\documentclass[a4paper,11pt]{article}
\usepackage[nointlimits]{amsmath}
\usepackage{amsfonts}
\usepackage{amssymb}
\usepackage{amsmath}
\usepackage{amsthm}
\usepackage{latexsym}
\usepackage{graphicx}
\usepackage{layout}
\usepackage[english]{babel}
\usepackage{color}

\numberwithin{equation}{section} 
        
\newtheorem{lemma1}     {Lemma}[section]
\newtheorem{teorema1}   [lemma1]{Theorem}
\newtheorem{prop1}      [lemma1]{Proposition}
\newtheorem{coroll1}    [lemma1]{Corollary}
\newtheorem{cong1}      [lemma1]{Conjecture}
\newtheorem{remark1}    [lemma1]{Remark}
\newtheorem{defin1}     [lemma1]{Definition}

\newenvironment{Lemma}[1][]
        {\begin{lemma1}[#1]\begin{samepage}}{\end{samepage}\end{lemma1}}
\newenvironment{Theorem}[1][]
        {\begin{teorema1}[#1]\begin{samepage}}{\end{samepage}\end{teorema1}}

\newenvironment{Remark}[1][]
        {\begin{remark1}[#1]\begin{samepage}}{\end{samepage}\end{remark1}}

\newcommand{\eps}     {\varepsilon}

\newcommand{\funtheta}     {\mathcal E}
\newcommand{\funhdue}     {v}
\newcommand{\growexp} {a}

\newcommand{\grad} {\nabla}
\newcommand{\Hess}       {\grad^2}
\newcommand{\HH}       {\mathcal{H}}

\newcommand{\jump} {J}

\newcommand{\cNN}      {\ensuremath{\mathbb N}}

\newcommand{\Om}       {\Omega}
\newcommand{\ou}       {(0,1)}

\newcommand{\PP}       {\mathcal{P}}

\newcommand{\um}     {u^-}

\newcommand{\up}     {u^+}

\newcommand{\gga}{\gamma}            
\newcommand{\gd}{\delta}

\newcommand{\gl}{\lambda}

\newcommand{\R}{\ensuremath{\mathbb R}}
\newcommand{\Z}{\ensuremath{\mathbb Z}}
\newcommand{\N}{\mathbb{N}}
\newcommand{\bS}{\mathbb S}

\newcommand{\cC}{\mathcal C}
\newcommand{\cH}{\mathcal{H}}

\newcommand{\cF}{\mathcal F}
\newcommand{\cA}{\mathcal A}

\newcommand{\cI}{\mathcal I}

\newcommand{\bbar}{\big{|}}

\newcommand{\dive}{\operatorname{div}}

\newcommand{\mesrest}{\text{\huge$\llcorner$}}

\newcommand{\sprod}[2]{\langle #1, #2 \rangle}

 \newcommand{\n}{n}
\newcommand{\Rn}{\R^\n}
\newcommand{\hausmoins}{\mathcal{H}^{\n-1}}
\newcommand{\Sn}{\bS^{\n-1}}
\newcommand{\Ldeu}{L^2(\Omega)}

\begin{document}

\title{The $\Gamma$-limit 
for singularly perturbed functionals of Perona-Malik type
in arbitrary dimension}
\author{Giovanni Bellettini
\thanks{Dipartimento di Matematica,
Universit\`a di Roma Tor Vergata,
via della Ricerca Scientifica 1, 00133 Roma, Italy,
and LNF-INFN,
via E. Fermi 40,  
00044 Frascati,  Italy.
E-mail: \texttt{belletti@mat.uniroma2.it}}
\and
Antonin Chambolle
\thanks{CMAP, Ecole Polytechnique, CNRS,
Palaiseau, France.\newline
E--mail: \texttt{antonin.chambolle@cmap.polytechnique.fr}}
\and
Michael Goldman
\thanks{
Max Planck Institute for Mathematics in the Sciences, Leipzig, Germany.  E-mail:  \texttt{goldman@mis.mpg.de}, partially funded by a Von Humboldt PostDoc fellowship}
}
\date{}

\maketitle

\begin{abstract}
In this paper we generalize to arbitrary dimensions
a one-dimensional equicoerciveness and 
$\Gamma$-convergence result 
for a second derivative 
perturbation of Perona-Malik type functionals.
Our proof relies on a new density result in the space of special functions of bounded variation with vanishing diffuse gradient part. This provides
a direction of investigation to derive approximation for functionals
with discontinuities penalized with a ``cohesive'' energy, that is,
whose cost depends on the actual opening of the discontinuity.
\end{abstract}

\section{Introduction}\label{secint}

We investigate in this paper a singular pertubation problem whose limit
is defined on piecewise constant functions, and corresponds essentially
to a subadditive penalization of the discontinuity.


More precisely, we consider for $\Omega$ a bounded open subset of $\R^n$ with Lipschitz boundary,
and $\nu \in (0,1]$, the functional $F_\nu : L^2(\Omega) \to [0,+\infty]$ defined by
\begin{equation}\label{ilmoneln}
F_\nu (u) := 
\left\{
\begin{array}{ll}
\displaystyle 
\frac{1}{2} \int_{\Omega} \left[ \nu^3 \vert \grad^2 
u\vert^2+\frac{1}{\nu \phi(1/\nu)}\phi(\vert \grad u\vert)
\right] ~dx 
& \mbox{if $u \in H^2(\Omega)$},
\\
\\
+ \infty & \mbox{elsewhere in $L^2(\Omega)$}.
\end{array}
\right.
\end{equation}
Here $\grad^2 u$ is the hessian of $u$ and, for a symmetric
real $(n\times n)$-matrix $M$, we set $\vert M\vert := \max
\{\langle M\xi, \xi\rangle: \xi \in \Rn, \vert \xi\vert \leq 1\}$.
Moreover the 
function $\phi:\R \to [0, +\infty)$ is continuous, 
even, nondecreasing, typically nonconvex and 
nonconcave with $\phi^{-1}(0)=\{0\}$, and 
 has sublinear growth at infinity,  in the sense that
\begin{equation}\label{eq:a}
\exists~ a \in [0,1) \textrm{ such that } \lim_{p\to +\infty} \frac{\phi(\gl p)}{\phi(p)}=\gl^a, \qquad \gl >0.
\end{equation}

In this paper we prove that the sequence $(F_\nu)$ 
is equicoercive and  
 $\Gamma$-converges \cite{Da:93,Br:02} 
as $\nu \to 0^+$ 
to the functional $\cF : L^2(\Omega) \to [0,+\infty]$ defined as 
follows:
\begin{equation}\label{thelimitn}
\cF (u) := 
\left\{
\begin{array}{ll}
\sigma_\growexp 
\displaystyle \int_{\jump_u} 
\vert \up - \um \vert^{\frac{2+\growexp}{4-\growexp}} ~ d\HH^{n-1}
& \mbox{if $u \in X(\Omega) \cap L^2(\Omega)$,}
\\
\\
+ \infty & \mbox{elsewhere in $L^2(\Omega)$},
\end{array}
\right.
\end{equation}
where $X(\Omega)$ consists of  special functions 
of bounded variation in $\Omega$ having gradient with no 
absolutely continuous part (see Section \ref{secnot}) and $\sigma_a>0$ is a constant depending only on $a$. 
Our results generalize an equicoerciveness and
$\Gamma$-convergence
result obtained in 
\cite{BelFus} in the one-dimensional case. 

The interest in these results are twofold. A first important application of such variational problems is in the numerical
analysis of fracture mechanics, and in particular the variational
approach to fracture growth popularized in the 90's by G.~Francfort
and J.~J.~Marigo, built upon the Mumford-Shah functional of image
processing~\cite{FrancfortMarigo,MumfordShah}. The original model
allows to approach the so-called ``Griffith'' model, where the
cost for opening a crack is proportional to its length or area surface.
Its success for predicting realistic fractures is impressive~\cite{BourdinFrancfortMarigo}, however its physical relevance is still a matter of discussion.
More physical models (known as ``cohesive'' and initially introduced by
Barenblatt), consider that for a fracture with a small opening,
the cost should rather be proportional (through some non-linear correspondence)
to the size of the discontinuity. Mathematically, the study of such
models is more tricky. Also, finding reasonable approximations
of such energies is a difficult problem  (still only partially solved,
if one really wants to consider linearized or non-linear elasticity
energies).
Some approaches need to consider a small non-cohesive term~\cite{AmbrosioLemenantRC,DalMasoIurlano}. 
A more physical term (with slope one when
the opening goes to zero) is obtained in~\cite{AlicandroBraidesShah},
which can approach quite general cohesive energies in the
scalar setting (and is, to our
knowledge, the best and most useful result in this direction so far).
It is built as a variant of Ambrosio and Tortorelli's approximation of the Mumford-Shah functional~\cite{AmTo}. 
The variant in~\cite{GobbinoMora}, built upon finite-differences
approximations, lies in between the previous results since it allows
to approach a partially cohesive term with an infinite slope for
infinitesimal openings.

This latter result is close to the model that
we present here, which is however built upon different ideas. Instead of
introducing a phase-field~\cite{AlicandroBraidesShah} or considering
finite differences, we penalize the variations of the gradient by
a higher order term. This is similar to the  two-wells problem
which has been studied in the celebrated paper~\cite{CFL:02}
of Conti, Fonseca, Leoni. Our energy
could be considered as a special case, where one well is
at $0$ and the other at infinity, just as~\cite{BoDuSe:96, OudSant} are to the
standard Modica-Mortola energy\footnote{
See also  \cite{AlMu, FoMa, contiben,Br:98}
for related problems.}. 

The second motivation is the relation of our result with 
the long-time behaviour of solutions
to the Perona-Malik equation \cite{PeMa:90}, 
obtained 
as the formal $L^2$-gradient flow of the functional
$$
{\rm PM}(u):=\frac{1}{2} \int_\Omega \log(1+|\nabla u|^2) \ dx,
$$
and
corresponding to the choice
\begin{equation}\label{perona}
\phi(x)=\log(1+\vert x\vert^2), \qquad x \in \Rn,
\end{equation}
and 
$a=0$ in \eqref{eq:a}. The Perona-Malik equation therefore reads as
\begin{equation}\label{peronamalikorigin}
u_t=\dive \left(\frac{\nabla u}{1+|\nabla u|^2} \right),
\end{equation}
and 
it is  ill-posed\footnote{
See for instance  
\cite{Ki:97,KaKu:98,gg,Go:03,BeNoPaTo06,BeNoPa:06,BeNoPaTo:07,Esedoglu}.}, 
due to the nonconvexity of $\phi$. 
In order to overcome the backward parabolic character
of \eqref{perona}, various regularizations 
have been suggested in the literature
(see \cite{BelFus} and  references therein);
in particular 
\cite{De:95}  
one can consider, for $\eps
\in (0,1]$, the 
functionals
\[
{\rm PM}_\eps(u):= \frac{1}{2} 
\int_\Omega \eps^2 |\nabla^2 u|^2+\log\left(1+|\nabla u|^2\right) \ dx,
\]
and take the corresponding gradient flow equations,
which amounts to add to the equation \eqref{peronamalikorigin} a 
fourth order term multiplied by $\eps^2$; see also \cite{slemrod}.
On the basis of the  numerical experiments \cite{BelFus} performed
 when $n=1$,  one observes 
various distinct time scales for the regularized equation; 
in particular, in a slow time scale,
the regularized solutions seem to converge to a 
piecewise constant function with the 
plateaus suitably evolving in the vertical direction.
The functionals $F_\nu$
are related 
to this 
slow time 
scale\footnote{We refer also to \cite{CoGo} for  the study of the slow time behaviour in a related semi-discrete approximation.}
setting
$\eps^2:=\nu^4 \log(1+\frac{1}{\nu^2})$
and taking $\frac{{\rm PM}_\eps}{\nu \log(1+\frac{1}{\nu^2})}$.
Therefore, one expects that
 the asymptotic limit as $\nu \to
0^+$ of the gradient flow of $F_\nu$ should shade some light 
on the behaviour of 
solutions to \eqref{peronamalikorigin}
for large times. This problem has been addressed in \cite{BelFus}
when $n=1$. The $n$-dimensional case seems much more difficult, 
and requires a preliminar study of the asymptotic limit of $F_\nu$
and of its variational properties, and this is the content of the present paper.
The most original part of this paper concerns the proof of 
the $\Gamma$-limsup inequality (Theorem \ref{pro:gsup}), 
which relies
on a new density result of simple plateaus functions in the space of special functions with 
bounded variation with vanishing diffuse gradient part,
see Lemma \ref{teoantonin}. 
Differently with respect
to the  
density theorem in \cite{CoTo:}, 
our result is valid without
 assumptions on the measure 
of the jump set of the limit functions. This allows us to get a full $\Gamma$-convergence result, differently from other related works \cite{AlGe:01, Mo:01},
see Remark \ref{rem:final}. 

\smallskip
 
The plan of the paper is the following. In Section \ref{secnot} we introduce
the notation and recall some definitions about $SBV$ 
functions and the slicing method. In Section \ref{secequigamma} we 
remind the results of \cite{BelFus} regarding 
the one-dimensional problem. 
In Section \ref{sec:thendim} we state and prove our density result 
and eventually in Section \ref{subsecndim} we 
prove the equicoerciveness and $\Gamma$-convergence theorems.

\section{Notation}\label{secnot}
In what follows $n \geq 1$, $\Omega \subset \Rn$ 
is a bounded open  Lipschitz set,
and  $\cA(\Omega)$ is the class of all open subsets of $\Om$. 
We denote by  $(e_1,\dots, e_\n)$ a fixed orthonormal basis of $\Rn$,
 by $\vert \cdot\vert$  (respectively $\langle \cdot, \cdot \rangle$) the euclidean  norm (respectively the euclidean scalar product)  in $\R^n$. We denote by $\mathcal H^{n-1}$ and  by $dx$, 
the $(n-1)$-dimensional Hausdorff measure  and the Lebesgue measure in $\Rn$.  The open ball of radius $\rho>0$ centered at $x \in \Omega$ is denoted by $B_\rho(x)$.
Throughout the paper, with a small abuse of 
language, we call sequence a family $(u_\nu)$ of functions 
labelled by a continuous parameter $\nu\in (0,1]$. A subsequence 
of $(u_\nu)$ is any sequence $(u_{\nu_k})$ such that $\nu_k
\to 0$ as $k \to +\infty$.

\subsection{$BV(\Omega)$ and $SBV(\Omega)$ functions}
$BV(\Om)$ is the space of functions 
$u \in L^1(\Om)$ 
 having as distributional derivative $Du$
a measure with finite total variation. For $u\in BV(\Om)$, we denote by $S_u$ the complement of the Lebesgue set of $u$. That is, $x \notin S_u$ if and only if 
$\lim_{\rho \to 0^+}\displaystyle \int_{B_\rho(x)} |u(y)-z| \ dy =0$
for some $z \in \R$.
 We say that $x$ is an approximate jump point of $u$ 
if there exist  $\xi \in \Sn$ 
and distinct $a, b \in \R$
 such that
\[ \lim_{\rho \to 0} \frac{1}{|B_\rho^+(x,\xi)|} \int_{B_\rho^+(x,\xi)} |u(y)-a| \ dy =0 \quad \textrm{ and } \quad  \lim_{\rho \to 0} \frac{1}{|B_\rho^-(x,\xi)|} \int_{B_\rho^-(x,\xi)} |u(y)-b| \ dy =0,\]
where $B_\rho^\pm(x,\xi):= \{ y \in B_\rho(x) : \pm \sprod{y-x}{\xi}>0 \}.$ Up to a permutation of $a$ and $b$ and a change of sign of $\xi$, this characterize the triplet $(a,b,\xi)$ which is then denoted by $(u^+,u^-, \nu_u)$. The set of approximated jump points is denoted by $J_u$. 
The following  theorem holds
\cite{AmFuPa:00}.
\begin{Theorem}
The set $S_u$ is countably $\hausmoins$-rectifiable and $\hausmoins(S_u\backslash J_u)=0$. Moreover $Du \mesrest J_u=(u^+-u^-) \nu_u \hausmoins \mesrest J_u$.
\end{Theorem}
We indicate by $Du= \grad u \ dx \ + \ D^s u$ the Radon-Nikodym decomposition of $Du$. 
Setting $D^c  u:= D^s u \mesrest (\Om \backslash S_u)$ we get the decomposition
\[Du= \grad u \ dx \ +\ (u^+-u^-)\nu_u \hausmoins \mesrest J_u\ +\  D^c u,\]
where $\mesrest$ denotes the  restriction.
When $n=1$ we use the symbol $u'$ in place of $\grad u$,
and $u(x^\pm)$ to indicate the right and left limits at $x$.
We let 
\[
\begin{aligned}
& SBV(\Om):=\{ u \in BV(\Om) :  D^c u=0\},
\\
& GSBV(\Om):=\{ u \in L^1(\Omega) : \max(-T,\min(T,u))\in SBV(\Om) \quad \forall T\in \R\},
\end{aligned}
\]
 and 
\[X(\Om):= \{ u \in SBV(\Om) : \nabla u=0\}.
\]

\subsection{Slicing}\label{sub:slicing}
In this section we recall the slicing method for functions with bounded variation. 
Let $\xi \in \Sn$ and let 
$$
\Pi_\xi:= \{ y \in \Rn : \sprod{y}{\xi}=0\}. 
$$
If $y \in \Pi_\xi$ and $E \subset \Rn$, we define  the one-dimensional
slice
$$
E_{\xi y}:=\{t \in \R : y+t \xi \in E\}.
$$
For $u: \Om \to \R$, we define $u_{\xi y} : \Om_{\xi y} \to \R$  as
\begin{equation}\label{eq:slice}
u_{\xi y}(t):=u(y+t \xi),  \qquad t \in \Om_{\xi y}.
\end{equation}
Functions in $GSBV(\Om)$ can be characterized 
by one-dimensional slices (see \cite[Thm. 4.1]{Br:98}).
\begin{Theorem}\label{teo:slice}
Let $u \in GSBV(\Om)$.
Then for all $\xi \in \Sn$ we have 
$$
u_{\xi y} \in GSBV(\Omega_{\xi y}) \qquad {\rm for}~ \hausmoins- {\rm a.e.}~ 
y \in \Pi_\xi. 
$$
Moreover for such $y$, we have
\begin{equation}\label{egalprim}
 u'_{\xi y }(t)= \sprod{\nabla u(y+t \xi)}{\xi} \quad \textrm{for a.e. } t \in \Om_{\xi y},
\end{equation}
\begin{equation}\label{jumpslice}
 J_{u_{\xi y}}=\{ t \in \R : y+t \xi \in J_u\},
\end{equation}
and
\begin{equation}\label{formjumpslice}
 u_{\xi y} (t^\pm)=u^\pm(y+t \xi) \quad \textrm{or} \quad u_{\xi y} (t^\pm)=u^\mp(y+t \xi),
\end{equation}
according to whether $\sprod{\nu_u}{\xi}>0$ or $\sprod{\nu_u}{\xi}<0$.
Finally, for every Borel function $g: \Omega \to \R$, 
\begin{equation}\label{integjumpslice}
 \int_{\Pi_\xi} \sum_{t \in J_{u_{\xi y}}} g_{\xi y} (t ) 
\ d\hausmoins (y)=\int_{J_u} g~ |\sprod{\nu_u}{\xi}| \ d\hausmoins.
\end{equation}
Conversely if $u \in L^1(\Om)$ and if for all $\xi \in \{e_1,\dots, e_n\}$ and  almost
every $y \in \Pi_\xi$ we have $u_{\xi y} \in SBV(\Om_{\xi y})$ and
\[\int_{\Pi_\xi} |Du_{\xi y}|(\Om_{\xi y})\ d\hausmoins(y) < +\infty,\]
then $u \in SBV(\Om)$. 
\end{Theorem}

\section{The one-dimensional case}\label{secequigamma}
In this section we briefly record the main results of \cite{BelFus},
obtained in dimension $n=1$,
that will be necessary in order to analyze the problem
in arbitrary dimension.
For 
$I\subset \R$ a bounded open interval 
we consider the functional
$$
F_\nu(u,I):= 
\begin{cases}
\displaystyle \frac{1}{2}
\int_I \left[\nu^3 (u'')^2+\frac{1}{\nu \phi(1/\nu)} \phi(u')\right] \ dx

& {\rm if}~ u \in H^2(I),
\\
\\
+\infty & {\rm if} ~ u \in L^1(I) \setminus H^2(I)
\end{cases}
$$
and 
\[
\cF(u,I):= 
\begin{cases}
\displaystyle \sigma_a \sum_{x\in J_u} |u(x^+)-u(x^-)|^{\frac{2+a}{4-a}} &
{\rm if}~ u \in X(I),
\\
+\infty & {\rm if}~ u \in L^1(I) \setminus X(I),
\end{cases}
\]
where $X(I)$ is defined, accordingly to the  $n$-dimensional case, as
$X(I) := \{u \in SBV(I) : u' =0\}$. 

Then the following results hold
\cite[Lemma 3.2 and Section 4]{BelFus}.
\begin{Lemma}\label{1Dcompact}
There exist $\overline \nu>0$, 
 a decreasing function $\omega: (0,\bar \nu) \to 
(0,+\infty)$ with $\displaystyle \lim_{\nu \to 0^+} \omega(s)=0$,
and a constant $C>0$ 
 such that for any  $u \in H^2(I)$ and $\nu \in (0,\overline \nu)$,
\begin{equation}\label{majvarunD}
 \int_I |u'|~ dx \le |I| \omega(\nu) +C F_\nu(u,I).
\end{equation}
 \end{Lemma}
\begin{Theorem}[{\bf Equicoerciveness}]\label{1Dcomp}
Let $(u_\nu)\subset H^2(I)$ be a sequence satisfying 
$$
\sup_{\nu \in (0,1]} F_\nu(u_\nu,I) < +\infty
$$
and such that $\displaystyle \int_I 
u_\nu~dx=0$ for any $\nu \in (0,1]$. Then there exist
a function $u \in X(I)$ and 
a subsequence of $(u_\nu)$ weakly*
converging 
 to $u$ 
in $BV(I)$.
\end{Theorem}

\begin{Theorem}[{\bf $\Gamma$-convergence}]\label{1Dgamma}
We have
\[\Gamma(L^1(I))-\lim_{\nu \to 0^+} F_\nu(\cdot,I) =\cF(\cdot,I),\]
and the same result holds true for the $\Gamma(L^2(I))$-limit.
\end{Theorem}

We also remind the construction of the 
recovery sequence in the proof of the $\Gamma-\limsup$
inequality in Theorem \ref{1Dgamma}. 
For any $\eta>0$ and $s\geq 0$, let
\[Y_\eta(s):= \left\{
\psi \in H^2((0,\eta)) : \psi(0)=0, \psi(\eta)=s, \psi'(0)=\psi'(\eta)=0\right\}.
\]
Define 
$$
e_a(s) := \frac{1}{2} 
\inf
\left\{
\int_{(0,\eta)} 
\Big[ (\psi'')^2  +  \vert \psi'\vert^a\Big]~ dx : 
\eta > 0, \psi \in Y_\eta(s)
\right\}.
$$
Then it turns out that 
$$
e_a(s) = \sigma_a s^{\frac{2+a}{4-a}},
$$
for a suitable constant $\sigma_a>0$.

Fix now $b>0$ and $s >0$.
Let $(\eta,\psi)\in (0,+\infty) \times Y_\eta(s)$ be such that
\[
\int_0^\eta \Big[
(\psi'')^2 + |\psi'|^a\Big] \ dx \le 
\sigma_a s^{~\frac{2+a}{4-a}} +b.\]
For 
a function of $X(I)$ of the form
$u= s~ 1_{[0,+\infty)}$, a recovery sequence $(u_\nu) \subset H^2(I)$ is given by
\[u_\nu(x):= \left\{\begin{array}{ll} 0 & x< 0,
\\
\psi(\frac{x}{\nu})& x \in [0, \eta \nu],
\\
s & x>\eta \nu,
\end{array}\right.\]
in the sense that
for every $\gd>0$ we have that
$(u_\nu)$ tends to $u$ in $L^1((-\gd,\gd))$ and
\[
\limsup_{\nu \to 0^+} 
F_\nu \left(u_\nu, (-\gd, \gd)\right)= \sigma_a s^{~\frac{2+a}{4-a}} + b=\cF(u,(-\gd,\gd)) +b,
\]
so that in particular
$\displaystyle \lim_{b \to 0^+}\limsup_{\nu \to 0^+} 
F_\nu \left(u_\nu, (-\gd, \gd)\right)=\cF(u,(-\gd,\gd))$.
A similar construction can be made for $s<0$.
\section{The $\n$-dimensional case: a density result}\label{sec:thendim}
Let $\theta: \R \to \R^+$ be an even lower semicontinuous function, 
subadditive and nondecreasing on $(0,+\infty)$, and  
continuous at $t=0$, such that 
$\theta(0) =0$, 
$\lim_{t \to 0} \frac{\theta(t)}{t} = +\infty$. Let us notice that the functions $s\to\vert s\vert^{\frac{2+a}{4-a}}$ satisfy these properties.

Define
$$
\PP(\Omega) := \{u \in GSBV(\Omega) : \nabla u= 0 ~ {\rm in}~ \Omega\}
$$
and for any  $A \in \cA(\Omega)$ set
\begin{equation}\label{thelimitant}
\funtheta(u,A) := 
\left\{
\begin{array}{ll}
\displaystyle \int_{\jump_u \cap A} \theta\left(u^+ - u^-\right)~d\HH^{n-1}
& \qquad \mbox{if $u 
\in \PP(A)$},
\\
\\
+ \infty & \qquad \mbox{elsewhere in $L^1(A)$}.
\end{array}
\right.
\end{equation}
We recall
that $\mathcal E(\cdot,A)$ 
is $L^1(A)$-lower semicontinuous \cite[Ex. 5.23]{AmFuPa:00}.

In this section we want to prove the 
following result, which is the main technical tool of this
paper, and 
will be used in the 
proof of Theorem \ref{pro:gsup}. 
\begin{Lemma}[{\bf Density}]\label{teoantonin}
Let $u \in \PP(\Omega)$ be such that $\funtheta(u,\Omega) < +\infty$. Then
there exists a sequence $(u_h) \subset \PP(\Omega)$ with the following
properties:
\begin{itemize}
\item[-] $\displaystyle
\lim_{h \to 0^+} u_h = u$ 
 in $L^1(\Omega)$; 
\item[-]  $\displaystyle \lim_{h \to 0^+}
\funtheta(u_h,\Omega) 
=\funtheta(u,\Omega)$;
\item[-] $\jump_{u_h}$ is contained in a finite 
union of facets of polytopes\footnote{A polytope is a set 
whose boundary consists  of a finite number of pieces of hyperplanes.} for any $h \in \cNN$. In 
particular
 for any $h \in \cNN$,
$$
\HH^{n-1}(\Omega \cap 
\overline \jump_{u_h}) = 
\HH^{n-1}(\Omega \cap \jump_{u_h}) \quad {\rm and}\quad
\HH^{n-1}(\overline \jump_{u_h}) < +\infty.
$$
\end{itemize}
\end{Lemma}
\begin{proof}
Let $\Omega' \supset \supset \Omega$ be an open set 
and let us denote still by $u \in \PP(\Omega')$ an extension of 
$u$ such that $\HH^{n-1}(\jump_u \cap \partial \Omega) 
=0$ (see \cite{Gi:84}).
We also fix a
representative of $u$ defined everywhere. 

We will need an auxiliary functional $\funtheta_0$ defined
as follows: for any open set $A \subseteq \Omega'$ 
\begin{equation}\label{thelimitantonin}
\funtheta_0 (v,A) := 
\left\{
\begin{array}{ll}
\displaystyle \int_{\jump_v \cap A} \Vert \nu_v \Vert_1 ~\theta(v^+ - 
v^-)~d\HH^{n-1}
& \qquad {\rm if}~ v 
\in \PP(A),
\\
\\
+ \infty & \qquad \mbox{elsewhere in $L^1(A)$},
\end{array}
\right.
\end{equation}
where, for $\xi = (\xi_1,\dots,\xi_n)\in
\Rn$, we let 
$$
\Vert \xi \Vert_1 := \vert \xi_1\vert + \dots +\vert \xi_n\vert.
$$
Again by \cite[Ex. 5.23]{AmFuPa:00}, 
the functional $\funtheta_0(\cdot,A)$ is $L^1(A)$-lower semicontinuous.

Observe that $\nu \in \mathbb S^{n-1}$ implies $\Vert \nu\Vert_1 \geq 1$,
with equality only if $\nu$ coincides with one of the vectors
of the orthonormal  basis
$(e_1,\dots,e_n)$ of $\R^n$. 
In particular
$$
\funtheta(u, \Omega) \leq \funtheta_0(u, \Omega).
$$
From our assumptions on $\theta$, we have that $\funtheta(u,\Omega)$ 
decreases under truncation. By the lower semicontinuity of 
$\funtheta(u,\Omega)$, this implies that setting  $u_T:=\max(-T,\min(T,u))$
for any $T>0$, 
then $\funtheta(u_T,\Omega)$ converges to $\funtheta(u,\Omega)$ 
as $T\uparrow +\infty$. Hence, with no loss of generality, we can 
assume that $u \in L^\infty(\Omega)$. 

We divide the proof into five steps. In the first
step we construct a suitable discrete approximation $u_\eps^y$ of $u$,
defined on points of a lattice.
For any 
$\delta >0$ set 
\[\Omega^\delta := \{x \in \R^n : 
{\rm dist}(x, \Omega) < \delta\}.\] 
{\it Step 1}. 
Let $\eps >0$ be such that $\Omega_\eps:= \Omega^{2 n^{1/2} \eps} \subset\subset 
\Omega'$
and, for any $y \in [0,1)^n$, set
$$
u_{\eps}^y(x) := u(x + \eps y), \qquad x \in \eps \Z^n \cap \Omega'.
$$
Define the discrete energy $D_\eps^y$  
of $u_\eps^y$ as
\begin{equation}\label{eq:defdiscrete}
D_\eps^y := \displaystyle 
\eps^{n-1} \sum_{i=1}^n 
~~ \sum_{x, x + \eps e_i \in 
\Omega_{\eps/2}
} 
\theta\big(u_\eps^y(x+\eps e_i) - u_\eps^y(x)\big).
\end{equation}
where $\Omega_{\eps/2}:= \Omega^{n^{1/2} \eps}$.
We claim that
\begin{equation}\label{primoclaim}
\int_{\ou^n} D_\eps^y ~dy \leq \funtheta_0(u, \Omega_\eps).
\end{equation}
To show the claim, we will follow some arguments 
in \cite{Go:98}, \cite{Ch:04} \cite{Ch:05} and \cite{ChGiPo:07}.
Let us first establish an inequality in one dimension: 
given a bounded open interval $I \subset \R$,
$\eps>0$ sufficiently small  and $v \in \PP(I)$, we have 
\begin{equation}\label{claimantonin}
\eps^{-1} \int_{I \cap (I - \eps)} 
\theta(v(t+\eps)- v(t))~dt \leq 
\sum_{t \in \jump_v} \theta(v^+(t) - v^-(t)).
\end{equation}
Indeed, for almost any $t \in I\cap (I -\eps)$, using the subadditivity
of $\theta$ it follows
$$
\theta(v(t+\eps)- v(t)) \leq \sum_{s \in \jump_v \cap (t, t+\eps)} 
\theta(v^+(s) - v^-(s)).
$$
Therefore, using Fubini's Theorem, 
\begin{align*}
\int_{I \cap (I - \eps)} 
\theta\left(v(t+\eps)- v(t)\right)~dt 
\leq &
\int_{I \cap (I - \eps)} 
\sum_{s \in \jump_v}
1_{(t, t+\eps)}(s) 
\theta\left(v^+(s) - v^-(s)\right)~dt 
\\
= & \sum_{s \in \jump_v}
\theta\left(v^+(s) - v^-(s)\right)
\int_{I \cap (I - \eps)} 
1_{(t, t+\eps)}(s) ~dt,
\end{align*}
which implies \eqref{claimantonin}. 

\smallskip

We can now prove claim \eqref{primoclaim}.
Since
$$
\int_{\ou^n} D_\eps^y ~dy = 
\eps^{n-1}
\sum_{i=1}^n 
~~\sum_{x, x + \eps e_i \in \Omega_{\eps/2}
} 
~\int_{\ou^n} 
\theta\Big(u(x+\eps y  + \eps e_i)) - u(x+\eps y)\Big)~dy,
$$
making the variable change $x' := \eps y + x$ (for fixed $x$)  
gives
\begin{align*}
\int_{\ou^n} D_\eps^y ~dy \leq & \eps^{-1} 
\sum_{i=1}^n 
\int_{\Omega_{\eps/2}
}
\theta\left(u(x'+\eps e_i) - u(x')\right) ~dx'
\\
= & 
\eps^{-1}
\sum_{i=1}^n
\int_{\Pi_{e_i}}
\int_{\left( \Omega_{\eps/2}\right)_{e_i \zeta}}
\theta\Big(u(\zeta + (t+\eps) e_i) - u(\zeta+ te_i)\Big) 
~ dt ~d\HH^{n-1}(\zeta).
\end{align*}
If, for any $i=1,\dots,n$ and $\mathcal H^{n-1}$-almost
every $\zeta \in \Pi_{e_i}$ we define as in 
Section \ref{sub:slicing}, 
the one-dimensional slice $u_{ e_i \zeta}$ as $u_{ e_i \zeta}(t)
:= u(\zeta + t e_i)$ for almost each $t \in \R$ such that $\zeta + t e_i 
\in 
\Omega'$, we get, 
$$
\int_{\ou^n} D_\eps^y ~dy \leq 
\eps^{-1}
\sum_{i=1}^n
\int_{\Pi_{e_i}}
\int_{\left(\Omega_{\eps/2}\right)_{e_i \zeta}}
\theta\Big(u_{e_i \zeta}(t+\eps) - u_{e_i \zeta}(t)\Big) 
~ dt ~d\HH^{n-1}(\zeta).
$$
Hence, using \eqref{claimantonin}, 
$$
\int_{\ou^n} D_\eps^y ~dy \leq 
\sum_{i=1}^n \int_{\Pi_{e_i}}
~\sum_{s \in \jump_{u_{e_i \zeta}} \cap
\left(\Omega_\eps\right)_{e_i \zeta}}
\theta\left(u^+_{e_i \zeta}(s) - u^-_{e_i \zeta}(s)\right)~d\HH^{n-1}(\zeta).
$$
Therefore, by \eqref{integjumpslice},
$$
\int_{\ou^n} D_\eps^y ~dy \leq
\sum_{i=1}^n 
\int_{\jump_u \cap \Omega_\eps}
\vert \langle \nu_u, e_i\rangle\vert
~\theta\left(u^+ - u^-\right)
~d\HH^{n-1} = \funtheta_0(u, \Omega_\eps),
$$
which is claim \eqref{primoclaim}. The proof of {\it Step 1} is concluded.

\medskip
In the next step we define
a suitable piecewise constant interpolation 
$\overline u_\eps^y$ of 
$u_\eps^y$ having the property that 
$\funtheta(\overline u_\eps^y,\Omega)$ and 
$\funtheta_0(\overline u_\eps^y,\Omega)$ coincide.

\smallskip

{\it Step 2}. 
Let $\eps$ and $y$ be as in {\it Step 1}. Define the 
function
$\overline u_\eps^y : \Omega \to \R$ as
\begin{equation}\label{eq:defue}
\overline u_\eps^y(z) 
:= \sum_{x \in \eps \Z^n}
u_\eps^y(x)~ 1_{C_\eps^y(x)}
(z), \qquad z \in \Omega, 
\end{equation}
where $C_\eps^y(x)$ denotes the open coordinate hypercube 
centered at $x + \eps y$ of side $\eps$, i.e., 
$$
C_\eps^y(x) := x + \eps y  + 
(-\eps/2, 
\eps/2)^n.
$$
Clearly $\overline u_y^\eps \in \PP(\Omega)$ and its jump set 
is contained in the union of the facets of the hypercubes 
$C_\eps^y(x)$. Let us prove that 
\begin{itemize}
\item[(i)] $\displaystyle \lim_{\eps \to 0} \int_{\ou^n} 
\Vert 
\overline u_\eps^y - u\Vert_{L^1(\Omega)}~dy =0$;
\item[(ii)] $\funtheta(\overline u_\eps^y, \Omega) 
= \funtheta_0(\overline u_\eps^y, \Omega) \leq D_\eps^y$.
\end{itemize}
To show (i) and (ii), 
we assume without loss of generality that 
all functions are extended to $0$ outside $\Omega'$. 
Making the variables change 
$$
(x,y) \in \eps \Z^n \times [0,1)^n \mapsto
y' := 
\eps y + x
$$
and 
then 
$$
y' \in \Rn \to \xi = z-y' \quad ({\rm for~ fixed~} z), 
$$
we obtain 
\begin{align*}
\int_{\ou^n}
\Vert 
\overline u_\eps^y - u\Vert_{L^1(\Omega)}~dy 
=&\int_{\ou^n} \int_{\R^n}
\big\vert
\sum_{x \in \eps \Z^n} u(\eps y+x) 1_{C_\eps^y(x)}(z) 
- u(z)\big\vert ~dz dy
\\
= &
\int_{\ou^n} \int_{\R^n}
\vert
\sum_{x \in \eps \Z^n} \big(u(\eps y+x) - u(z)\big)~ 1_{C_\eps^y(x)}(z) \vert ~dz 
dy
\\
\leq& 
\int_{\ou^n} \sum_{x \in \eps \Z^n} 
~ \int_{\R^n}
\vert u(\eps y+x) - u(z)\vert ~ 1_{C_\eps^y(x)}(z)  
~dz dy
\\
= &
\eps^{-n}
\int_{\R^n}  \int_{\R^n} 
\vert u(y') - u(z)\vert 
~1_{(-\frac{\eps}{2},\frac{\eps}{2})^n}(z-y')  
~dz dy'
\\
=&
\eps^{-n}
\int_{\R^n} \int_{\R^n} 
\vert u(z-\xi) - u(z)\vert 
~1_{(-\frac{\eps}{2},\frac{\eps}{2})^n}(\xi)  
~dz d\xi
\\
= &
\eps^{-n}
\int_{(-\frac{\eps}{2},\frac{\eps}{2})^n}
\Vert u(\cdot - \xi) - u(\cdot)\Vert_{L^1(\R^n)} 
~d\xi
\\
= &
\int_{(-\frac{1}{2},\frac{1}{2})^n}
\Vert u(\cdot - \eps\xi) - u(\cdot)\Vert_{L^1(\R^n)} 
~d\xi,
\end{align*}
which tends to zero as $\eps \to 0$ by the dominated convergence theorem. 
This proves (i).

The equality in (ii) follows from the fact that 
$$
\Vert 
\nu_{
\overline u^y_\eps
}(x) \Vert_1^{} 
=1 \quad {\rm for~} \HH^{n-1}-{\rm a.e.}~x \in \jump_{\overline u_\eps^y},
$$
since the jump of $\overline u^y_\eps$ is contained in facets of coordinate
hypercubes. 

The inequality in (ii) follows from 
the definition of $\funtheta_0$ (see \eqref{thelimitantonin}) and 
definitions \eqref{eq:defdiscrete} of $D_\eps^y$ and
 \eqref{eq:defue} of $\overline u_\eps^y$. 

\medskip

{}From (i) we deduce that there exists a set $A \subset 
\ou^n$ with $\vert A\vert =1$ and a subsequence $(\eps_k) \subset (0,1)$ 
converging to zero such that 
$\overline u_{\eps_k}^y \to u$ in $L^1(\Omega)$ as $k \to \infty$ for any 
$y \in A$. 

{}From {\it Step 1} and Fatou's lemma we have
$$
\int_{\ou^n}
\liminf_{k \to +\infty} D_{\eps_k}^y~dy 
\leq \funtheta_0(u, \Omega).
$$
Hence there exists a point $\overline y \in A$ and a further subsequence
(still denoted by $(\eps_k)$) such that in addition 
there exists the $\lim_{k \to +\infty} D_{\eps_k}^{\overline y}$
and 
$$
\lim_{k \to +\infty} D_{\eps_k}^{\overline y} \leq \funtheta_0(u, \Omega).
$$
Summarizing, we have shown that given $u\in \PP(\Omega) \cap 
L^\infty(\Omega)$, if we set
$$
u_k := u^{\overline y}_{\eps_k}, \qquad k \in \N,
$$
the sequence $(u_k) \subset \PP(\Omega)$ has the following properties:
\begin{itemize}
\item[(a)] $\Vert u_k\Vert_{L^\infty(\Omega)} \leq
\Vert u\Vert_{L^\infty(\Omega)}$ for any $k \in \cNN$,
\item[(b)] $\displaystyle \lim_{k \to \infty} \Vert u_k - u\Vert_{L^1(\Omega)} =0$, 
\item[(c)] $\displaystyle \limsup_{k\to 
\infty} \funtheta_0(u_k, \Omega) \leq \funtheta_0(u,\Omega)$, 
\item[(d)] $\funtheta_0(u_k, \Omega) = \funtheta(u_k,\Omega)$ 
and $\jump_{u_k}$ 
is contained in a finite union of facets of coordinate hypercubes for any 
$k$.
\end{itemize}
To prove the theorem, it remains to show that we can replace $\funtheta_0(u, \Omega)$ with 
$\funtheta(u,\Omega)$ on the right hand side of the inequality in assertion (c). 

\medskip

{\it 
Step 3}. Let $A \subset \R^n$ be a bounded open set with 
Lipschitz boundary and let 
$v \in \PP(A) \cap  L^\infty(A)$ be with $\funtheta(v,A) <+\infty$. Let $(v_k)
 \subset 
\PP(A)$ be a sequence converging to $v$ in 
$L^1(A)$ such that $\Vert v_k\Vert_{L^\infty(A)} \leq \Vert 
v\Vert_{L^\infty(A)}$ for any $k \in \cNN$ and $\lim_{k \to \infty} 
\funtheta_0(v_k,A) = 
\funtheta_0(v,A)$. Then 
\begin{equation}\label{convtrace}
\lim_{k \to \infty}
\int_{\partial A} \vert T v_k - Tv\vert ~d\HH^{n-1} =0,
\end{equation}
where $Tv_k$ (resp. $Tv$)
denotes the trace of $v_k$ (resp. of $v$) on $\partial A$.

This follows from the fact that,
if $v \in \PP(A) \cap L^\infty(A)$, then 
 there exists a constant $c>0$ depending
only $\theta$ and $\Vert v\Vert_{L^\infty(A)}$ such that 
$\int_B \vert Dv\vert \leq c \funtheta_0(v,B)$ for  any open set $B\subseteq A$. 
A careful analysis 
of \cite[Theorem 3.88]{AmFuPa:00} (see also 
\cite[Section 5.3]{EvGa:}),
the $L^1(A)$-lower semicontinuity of 
$v \to \funtheta_0(v,B)$ for any open set $B \subseteq A$ and the fact that
$\funtheta_0(v,\cdot)$ is a regular measure, imply
\eqref{convtrace}.

\medskip

{\it Step 4}. Let $\mu$ be a Radon measure on $\Omega$ with values in 
$\R^d$, $d \geq 1$, with total variation measure
$\vert \mu\vert$ and $\vert \mu\vert(\Omega)< +\infty$. For any $h \in \cNN$,  let us consider a family
$\{B_k^h\}$, 
where $B_k^h$ are nonempty Borel subsets of $\Omega$ with $B_{k_1}^{h} 
\cap 
B_{k_2}^{h} = \emptyset$ when $k_1, k_2 \in \cNN$ and $k_1 \neq k_2$, 
and 
\begin{equation}\label{etauno}
\displaystyle \vert \mu\vert\Big(\Omega \setminus \bigcup_k B_k^h
\Big)=0.
\end{equation}
If in addition 
\begin{equation}\label{etadue}
\lim_{h\to +\infty} ~ ~
\sup_{k \in \cNN} ~{\rm diam}(B_k^h)=0,
\end{equation}
then 
\begin{equation}\label{vartota}
\vert \mu\vert (\Omega) = 
\lim_{h \to +\infty} 
\sum_k \vert \mu(B_k^h)\vert. 
\end{equation}
To prove \eqref{vartota} it is enough to show that for any $\varphi \in 
(\cC_c(\Omega))^d$ with $\Vert 
\varphi \Vert_{L^\infty(\Omega)} \leq 1$, we 
have 
$$
\int_\Omega \langle \varphi, \mu \rangle 
\leq \lim_{h \to + \infty} \sum_k \vert \mu(B_k^h)\vert.
$$
For any $h,k
\in \cNN$ pick a point 
$x_k^h \in B_k^h$, and define 
$$
\varphi_h(x) := \sum_k \varphi(x_k^h) ~1_{B_k^h}(x), \qquad x \in \Omega.
$$
Using the   
uniform continuity of $\varphi$ 
and assumption \eqref{etadue}, we have $\lim_{h \to +\infty} 
\Vert \varphi-\varphi_h\Vert_{L^\infty(\Omega,\mu)} =0$. 
Therefore
\begin{equation}
\int_\Omega \langle \varphi, \mu \rangle = 
\int_\Omega \langle \varphi_h, \mu \rangle 
+ \int_\Omega \langle (\varphi - \varphi_h),
\mu\rangle \leq 
\sum_k \vert \mu(B_k^h)\vert + o(1),
\end{equation}
and \eqref{vartota} follows.
\medskip

In the next (final) step we show that
it is possible to 
replace $\funtheta_0(u, \Omega)$ with $\funtheta(u,\Omega)$ in (c),
 by suitably
modifying the sequence $(u_k)$.

\smallskip

{\it Step 5}.  Let $\delta>0$ and $y \in (0,\delta)^n$. 
For each $z \in \delta \Z^n$ 
let us denote by $C_\delta(z)$ the open hypercube $y + z + 
(0,\delta)^n$. Possibly modifying the choice 
of the origin $y$, we can assume that
\begin{equation}\label{zeromeasbdry}
\displaystyle\vert Du\vert\Big(\Omega \setminus \bigcup_{z \in \gd\Z^n} 
C_\delta(z)\Big) =0.
\end{equation}
We will apply the construction described in {\it Steps 1,2} to
each of the sets $\Om \cap C_\delta(z)$ and then glue together the
sequences thus obtained.
Notice that it 
may happen that the set
$\Omega \cap C_\delta(z)$ is not Lipschitz (if non empty); 
however we already know that $u$ is defined on $\Omega'$ which 
strictly contains $\Omega$.  
Let us define the vector-valued Radon measure $\mu$ as
$$
\mu (B) := \int_{B\cap \jump_u}\theta\left(u^+ - u^-\right) \nu_u~ d\HH^{n-1},\qquad 
B {\rm ~Borel~ set~} \subseteq	\Omega.
$$
Observe that
 $\vert \mu\vert(B) = 
\int_{B \cap \jump_u} \theta\left(u^+-u^-\right)~d\HH^{n-1}$.

With any hypercube $C_\delta(z)$ we 
associate a unit vector $\xi_z \in S^{n-1}$ such that 
\begin{equation}\label{mispi}
\langle \xi_z,  
\mu\big(\Omega \cap C_\delta(z)\big)\rangle = \vert \mu\big(\Omega \cap 
C_\delta(z)\big)\vert.
\end{equation}
If we consider an orthonormal basis of $\R^n$ having $\xi_z$ as the 
first vector and if we let $\Vert\cdot\Vert_{1,z}$ be the $\Vert\cdot\Vert_1$ 
norm in this basis, from {\it Steps 1} and {\it 2}
 we can find a sequence $(u_k^\delta) \subset
L^1(\Omega \cap C_\delta(z))$ converging to $u$ 
in $L^1(\Omega \cap C_\delta(z))$ and satisfying properties 
(a)-(d) with $\Omega \cap C_\delta(z)$ 
in place of $\Omega$,
and such that 

\begin{align*}
\limsup_{k \to \infty} \funtheta(u_k^\delta, \Omega \cap C_\delta(z))
 \leq &\int_{\jump_u \cap C_\delta(z)} 
\Vert\nu_u\Vert_{1,z} ~d\vert \mu\vert
\\
\leq& \int_{\jump_u \cap C_\delta(z)} \vert \langle \xi_z,\nu_u\rangle\vert ~d\vert \mu\vert + c\int_{\jump_u \cap C_\delta(z)} \sqrt{1 - \langle \xi_z, \nu_u\rangle^2}~d\vert \mu\vert
\\
=& \mathcal E(u, \Omega \cap C_\delta(z))
+ c\int_{\jump_u \cap C_\delta(z)} \sqrt{1 - \langle \xi_z, \nu_u\rangle^2}~d\vert \mu\vert,
\end{align*}
where $c>0$ is a constant controlling the euclidean norm
of a vector with its $\Vert \cdot\Vert_1$-norm.
Using the same construction for all $z$ such that $\jump_u \cap 
C_\delta(z) \neq \emptyset$, 
and gluing together the functions obtained in each hypercube, 
we construct $u_k \in L^1(\Omega)$ such that 
\begin{equation}\label{primastima}
\funtheta(u_k, \Omega) = 
\sum_{z \in \delta \Z^n} 
\funtheta(u_k, \Omega \cap C_\delta(z))
+ \int_{\displaystyle \Omega \setminus 
\bigcup^{}_{z\in \delta \mathbb Z^n} C_\delta(z)}
~ \theta\left(u_k^+-u_k^-\right) ~d\HH^{n-1}.
\end{equation}
We notice that the jump set of the $u_k$'s consists of the union
of the facets of the lattice and of the jump set inside each
hypercubes. Recalling \eqref{zeromeasbdry} and {\it Step 3}, possibly
passing to a (not relabelled) subsequence, we can assume that 
the traces of $u_k$ on both sides of each facet
of a hypercube $C_\delta(z)$ converge $\vert \mu\vert$-almost
everywhere to the same 
limit
(which is the trace of $u$) as $k \to +\infty$, hence the last
term on the right hand side of \eqref{primastima} tends to zero
by the continuity of $\theta$ at $0$. 
Therefore there exists
$\overline k \in \cNN$ such that, if $k > \overline k$, then  
$\Vert u_k - u\Vert_{L^1(\Omega)} 
< \delta$ and 
\begin{equation}\label{dodici}
\begin{aligned}
\funtheta(u_k,\Omega) \leq &
\sum_{z \in \delta \Z^n} \left[
\int_{\jump_{u} \cap C_\delta(z)} \vert \langle \xi_z,
\nu_{u}\rangle\vert ~d\vert \mu\vert + 
c \int_{\jump_{u} \cap C_\delta(z)} \sqrt{1 - \langle \xi_z, 
\nu_{u}\rangle^2}
~d\vert \mu\vert  \right] + \delta 
\\
\leq & \funtheta(u, \Omega) + 
c \sum_{z \in \delta \Z^n} 
\int_{\jump_{u} \cap C_\delta(z)} \sqrt{1 - \langle \xi_z, 
\nu_{u}\rangle^2}
~d\vert \mu\vert + \delta,
\end{aligned}
\end{equation}
where in the last inequality
we made use of \eqref{mispi}.
The proof is therefore achieved if we show that the last sum on the right
hand side of \eqref{dodici} can be made small for $\delta>0$ sufficiently small.
We have
\begin{align}\label{tredici}
& \sum_{z\in \delta \Z^n} \int_{\jump_{u} \cap C_\delta(z)} 
\sqrt{1 - \langle \xi_z, 
\nu_{u}\rangle^2}
 \ d\vert \mu\vert 
\nonumber
\\
= &
\sum_{z\in \delta \Z^n} \int_{\jump_{u} \cap C_\delta(z)} 
\sqrt{(1 - \langle\nu_u,\xi_z \rangle)(1+\langle \nu_u,\xi_z\rangle)}  \
d \vert \mu\vert 
\\
\leq & \left( \sum_{z \in \delta \Z^n}
\int_{\jump_u \cap C_\delta(z)}
\Big(1 - \langle \nu_u, \xi_z\rangle\Big)~ \ d \vert \mu \vert 
\right)^{1/2}
\sqrt{2 \vert \mu\vert(\Omega)}.
\nonumber
\end{align}
Observing that from \eqref{mispi} 
$$
\int_{\jump_u \cap C_\delta(z)}
\Big(1 - \langle \nu_u, \xi_z\rangle\Big)~ d \vert \mu \vert \leq
\vert\mu \vert(\Omega) - \sum_{z \in \delta \Z^n} \vert 
\mu(\Omega \cap C_\delta(z))\vert,
$$
by {\it Step 4} we see that the right hand side of \eqref{tredici}
can be made arbitrarily small if $\delta>0$ is small enough. This concludes
the  proof.
\end{proof}

\section{Equicoerciveness and 
$\Gamma$-convergence 
}\label{subsecndim}

Using a slicing argument and Theorem \ref{teoantonin},
 the results of Section 
\ref{secequigamma} 
can be generalized in $n$-dimensions.

We remind that $F_\nu : \Ldeu \to [0,+\infty]$ 
and $\cF: \Ldeu \to [0,+\infty]$
are defined in \eqref{ilmoneln} and \eqref{thelimitn}, respectively.
If $A \in \mathcal A(\Omega)$ and $u \in H^2(A)$
we define the 
localized functional $F_\nu(u,A)$ as
$$
F_\nu(u,A) := \displaystyle 
\frac{1}{2} \int_{A} \left[ \nu^3 \vert \grad^2 
u\vert^2+\frac{1}{\nu \phi(1/\nu)}\phi(\vert \grad u\vert)
\right] ~dx.
$$
Similarly,
if $u \in X(A) \cap L^2(A)$, we set
$$
\cF(u,A) := \sigma_\growexp 
\displaystyle \int_{A \cap \jump_u} 
\vert \up - \um \vert^{\frac{2+a}{4-a}} ~ d\HH^{n-1}.
$$
We will  need the following lemma on the supremum of a family of measures (see \cite{Br:98}).
\begin{lemma1}\label{lemmeasur}
Let $\mu : \cA(\Om)\to [0,+\infty)$ be 
a superadditive set function
and  let $\gl$ be a positive measure on $\Om$.
For any $i \in \N$ 
let $\psi_i$ be a Borel function on $\Omega$ such that $\mu(A)\ge \int_A \psi_i ~d\gl$ 
for all $A\in \cA(\Om)$.
Then 
\[\mu(A) \ge \int_A \psi \ d\gl \qquad \forall A \in \cA(\Om)\]
where $\psi:=\displaystyle \sup_{i\in \N} \psi_i$. 
\end{lemma1}

We start our $n$-dimensional analysis 
with the following two results, which are independent of Lemma \ref{teoantonin}.

\begin{Theorem}[{\bf Equicoerciveness}]\label{equicoer}
 Let $(u_{\nu}) \subset H^2(\Omega)$ be a sequence satisfying
\begin{equation}\label{eq:bddene}
\sup_\nu F_{\nu}(u_{\nu}, \Omega) <+\infty
\end{equation}
and such that $\displaystyle \int_\Omega u_\nu ~dx =0$.
Then
there exist a function 
 $u \in X(\Omega)\cap \Ldeu$ and a 
subsequence of $(u_\nu)$ 
weakly$^*$ converging to $u$ in $BV(\Omega)$.
\end{Theorem}

\begin{proof}
 By Lemma \ref{1Dcompact}, 
there exist $\bar \nu\in (0,1]$, a decreasing function 
$\omega: (0,\bar \nu)\to (0,+\infty)$ with 
$\lim_{\nu \to 0^+} \omega(s)=0$ and a constant 
$C>0$ such that 
for every interval $I\subset \R$, every $v \in H^2(I)$ and
$\nu \in (0,\overline \nu)$,
inequality \eqref{majvarunD} holds, namely, for $\nu \in (0,\overline\nu)$,
\begin{equation*}
 \int_I |v'|~ dx \le |I|~ \omega(\nu) +C F_\nu(v,I).
\end{equation*}
Therefore, 
if $ \funhdue
\in H^{2}(\Omega)$ and $\xi \in \Sn$, by \eqref{minornu},
\begin{equation}\label{majvarxi}
\begin{aligned}
 \int_{\Omega} |\sprod{\nabla \funhdue}{\xi}| dx & =\int_{\Pi_\xi} 
\int_{\Omega_{\xi y}} |\sprod{ \nabla \funhdue(y+t \xi) }{\xi}| 
~dt \ d \hausmoins(y)
\\
&=\int_{\Pi_\xi} \int_{\Omega_{\xi y}} |\funhdue'_{\xi y}| ~dt 
\ d\hausmoins(y)
\\
&\le \int_{\Pi_\xi} \Big(
|\Omega_{\xi y}|~ \omega (\nu)+ CF_\nu(\funhdue_{\xi y}, \Omega_{\xi y}) \Big)
~d\hausmoins(y)
\\
&\le |\Omega|~ \omega(\nu) +C F_\nu(\funhdue,\Omega).
\end{aligned}
\end{equation}
Letting $\xi = e_i$ 
and summing 
 \eqref{majvarxi} over $i\in \{1,\dots,n\}$, we deduce
using \eqref{eq:bddene}, that there exists a positive
constant $\kappa$ such that 
\begin{equation}\label{majvarND}
 \sup_\nu \int_{\Omega} |\nabla \funhdue_\nu| \ dx 
\le \kappa \sup_\nu 
\Big(|\Omega|~ \omega(\nu) + F_\nu(\funhdue_\nu,\Omega)\Big) < +\infty.
\end{equation}
Thus if $\displaystyle \int_\Omega u_\nu ~dx=0$ 
then there exist
a function $u \in BV(\Omega)$
and a subsequence of $(u_\nu)$ 
weakly$^*$ converging to $u$ in $BV(\Omega)$.
\end{proof}

\begin{Theorem}[{\bf $\Gamma$-liminf}]\label{ah}
We have 
\begin{equation}\label{eq:liminfn}
\cF(\cdot,A)\le \Gamma(\Ldeu)-\lim_{\nu\to 0^+}
 F_\nu(\cdot,A), \qquad A\in \cA(\Om).
\end{equation} 
\end{Theorem}
\smallskip
\begin{proof}
The proof follows closely \cite[Prop. 3.4]{AlGe:01}. 
Fix $\xi \in \Sn$, $u \in H^2(\Omega)$,
$A \in \cA(\Omega)$
and let as usual $u_{\xi y} : A_{\xi y} \to \R$ be the slice 
 defined as in Section \ref{sub:slicing}, i.e.,
\[u_{\xi y}(t):=u(y+t \xi), \qquad t \in A_{\xi y}.\]
Since for almost every $t\in A_{\xi y}$ and $\mathcal
H^{n-1}$-almost every $y$  
\[
u_{\xi y}'(t)=\sprod{\nabla u(y+t \xi)}{ \xi} 
\qquad \textrm{ and } \qquad u_{\xi y}''=\sprod{\Hess u(y+t \xi)\xi}{\xi},\]
by Fubini's Theorem and being $\phi$ nondecreasing,
\begin{equation}\label{minornu}
\begin{aligned}
& F_\nu(u,A) \\
=& \int_{\Pi_\xi}\left[ \frac{1}{2} 
\int_{A_{\xi y}} \nu^3 |\nabla^2 u(y+t \xi)|^2
+\frac{1}{\nu \phi(1/\nu)} \phi(|\nabla u(y+t \xi)|)  
\ dt \ \right]d\hausmoins(y) 
\\
 \ge &\int_{\Pi_\xi} \left[\frac{1}{2} \int_{A_{\xi y}} \nu^3 
(\sprod{\nabla^2 u(y+t \xi) \xi}{\xi})^2
+\frac{1}{\nu \phi(1/\nu)} 
\phi(\sprod{\nabla u(y+t \xi)\xi}{\xi})  \ dt \right]\ d\hausmoins(y)
\\
=&\int_{\Pi_\xi} F_\nu(u_{\xi y}, A_{\xi y}) \ d \hausmoins(y).
\end{aligned}
\end{equation}

Let now $(u_\nu)\subset H^2(A)$ be a sequence
converging to $u$ in $L^2(A)$ as $\nu \to 0^+$. 
By
Fubini's theorem and Fatou's lemma
we have  $(u_\nu)_{\xi y}\to u_{\xi y}$ in $L^2(A_{\xi y})$ for 
$\mathcal H^{n-1}$-almost
every  $y \in \Pi_\xi$. Hence by \eqref{minornu} and Theorem \ref{1Dgamma},
\[
\sigma_a 
\sum_{ x \in J_{u_{\xi y}} \cap A_{\xi y}} 
|u_{\xi y}(x^+)-u_{\xi y} (x^-)|^{\frac{2+a}{4-a}}
\leq
\liminf_{\nu \to 0^+} 
F_\nu\left( (u_\nu)_{\xi y}, A_{\xi y}\right). 
\]
Thus, applying once more  Fatou's lemma,
\begin{equation}\label{liminfxi}
\int_{\Pi_\xi}\left( 
\sigma_a \sum_{ x \in J_{u_{\xi y}} \cap A_{\xi y}} |u_{\xi y}(x^+)
-u_{\xi y}(x^-)|^{\frac{2+a}{4-a}}\right) d \hausmoins(y) \leq
\liminf_{\nu \to 0^+} 
F_\nu(u_\nu,A).
\end{equation}
Let us 
first show that 
\begin{equation}\label{first}
u \in X(\Omega). 
\end{equation}
Let $T>0$ and set $u_T:=\max(-T,\min(T,u))$. By the results of Section 
\ref{secequigamma}
we know that $(u_T)_{\xi y } \in X(A_{\xi y})$ and 
thus by \eqref{liminfxi} and \eqref{eq:bddene},
\[\int_{\Pi_\xi} | D(u_T)_{\xi y}|(A_{\xi y}) 
\ d \hausmoins(y) <+\infty, \qquad 
\xi \in \Sn, \ \mathcal H^{n-1}-
{\rm a.e.} ~ y \in \Pi_\xi. 
\]
Therefore, by Theorem \ref{teo:slice} 
it follows that 
 $u \in GSBV(A)$. By \eqref{majvarND} we also know that $u \in BV(A)$ 
and thus $u \in SBV(A)$. Now for every $\xi \in \Sn$, 
\begin{align*}\int_A |\sprod{\nabla u}{\xi}| dx&= \int_{\Pi_\xi} \int_{A_{\xi y}} |\sprod{\nabla u(y+t\xi)}{ \xi}| ~dt \ d\hausmoins(y)\\
 						&=\int_{\Pi_\xi} \int_{A_{\xi y}} |u'_{\xi y} (t)| ~dt \ d\hausmoins(y)\\
						&=0.
\end{align*}
Hence $\nabla u =0$ almost everywhere, so that $u \in X(A)$. This concludes the proof of \eqref{first}. 

\smallskip

By Theorem \ref{teo:slice} (a), we get 
\begin{equation}\label{gammaliminfxi}
 \liminf_{\nu \to 0^+} F_{\nu} (u_\nu,A) \ge \sigma_a \int_{A\cap J_u} 
|u(x^+)-u(x^-)|^{\frac{2+a}{4-a}} |\sprod{\nu_u}{\xi}| \ d\hausmoins(y).
\end{equation}
 We now let $\gga$ be the 
superadditive increasing set function defined on $\cA(\Omega)$ by
\[\gga(A):=\Gamma-\liminf_{\nu \to 0^+} F_\nu (u,A)\]
and we let $\lambda$ be the Radon measure defined as 
\[\gl := | u(x^+)-u(x^-)|^{\frac{2+a}{4-a}}~ \hausmoins {\mesrest J_u}.\]
Fix a sequence $(\xi_i)_{i\in \N}$ dense in $\Sn$.
By \eqref{gammaliminfxi} we have
\[\gga(A) \ge \int_A \psi_i \ d\gl, \quad i \in \N,
\]
where 
\[
\psi_i(x) := 
\begin{cases} 
|\sprod{\nu_u(x)}{\xi_i}| & \textrm{if } x \in J_u\\[8pt]
                                       0 \qquad & \textrm{if } x \in A \backslash J_u.
                                                                              \end{cases} \]
Hence by Lemma \ref{lemmeasur}, 
\begin{align*}
 \Gamma-\liminf_{\nu \to 0^+} F_\nu (u,A)\ge \int_A \sup_i \psi_i(x) \ d\gl
					=\cF(u,A).
\end{align*}

\end{proof}

We next pass to the 
proof of the $\Gamma$-limsup inequality. 
We start by studying the particular case  $u= s 1_E$,
 adapting the proof from \cite[Prop 3.5]{AlGe:01}.

\begin{prop1}\label{gammasupE}
 Let $\Omega' \supset 
\supset \Omega$ be a bounded open set and let $E$ be such that $E=E' \cap \Omega$ where $E'$ is a set of finite perimeter in $\Omega'$ such that $\partial E' \cap \Omega'$ is a smooth hypersurface. Then
\begin{equation}\label{gammasupz}
 \Gamma-\limsup_{\nu \to 0^+} 
F_\nu( s 1_E, A) \le \sigma_a |s|^\frac{2+a}{4-a} 
~\hausmoins(A \cap \partial E),
\end{equation}
for all $s \in \R$ and $A \in \cA(\Omega)$.
\end{prop1}

\begin{proof}
 Let 
\[d(x) := {\rm dist}(x, \Rn \setminus  E')-{\rm dist}(x,E'), \qquad x\in \Rn,
\]
be the signed distance function from $\partial E'$ positive in 
$E'$.
Take  $\gd>0$ such that $d \in 
\cC^\infty(\overline{V_\gd \cap \Omega})$ where $V_\gd:= \{ x \in \Rn \ : 
\ |d(x)|<\gd\}$ (see for instance \cite{AmDan:00}),
and recall that $\vert \nabla d\vert^2=1$ in $V_\delta \cap \Omega$.
 Let $b>0$ and $(u_\nu)$ be the one-dimensional recovery sequence given in 
Section \ref{secequigamma} converging to 
$s 1_{\{t>0\}}$ in $L^2((-\delta,\delta))$, 
so that 
\begin{equation}\label{gammalimsupund}
\lim_{\nu \to 0^+} 
F_\nu\left(u_\nu, (-\gd,\gd)\right)=\sigma_a |s|^\frac{2+a}{4-a}+b.
\end{equation}
We next define
\[\widetilde u_\nu(x):= u_\nu(d(x)), \qquad 
x \in \Omega.
\]
Then $\widetilde u_\nu \in H^2(\Omega)$ 
and $\widetilde u_\nu \to s 1_E$ in $\Ldeu$
as $\nu \to 0^+$. Using the coarea formula, 
we find
\begin{eqnarray*}
 F_\nu(\widetilde u_\nu,A)
&= \displaystyle \frac{1}{2} &\int_{A 
\cap V_{\eta \nu}} \left[ 
\nu^3|u_\nu''(d) (\nabla d \otimes \nabla d) +  u_\nu'(d) \nabla^2 d|^2 +\frac{1}{\nu \phi(1/\nu)} \phi(|u_\nu'(d) \nabla d|) \right] \ dx
\\[8pt]
&\le 
\displaystyle\frac{1}{2}&
\int_0^{\eta \nu} \int_{\{x\in A:d(x)=t\}}
 \left[\nu^3(u_\nu''(t))^2  + 
2\nu^3|u_\nu''(t)||u_\nu'(t)|\bbar
 +\nu^3(u_\nu'(t))^2 \vert\vert \nabla^2 d \vert\vert_\infty^2
\right.\\
 &  & \left.  +\frac{1}{\nu \phi(1/\nu)} \phi(|u_\nu'(t)|) \right] d \hausmoins
~ dt\\[8pt]
&=\displaystyle\frac{1}{2} &
\int_{0}^{\eta \nu} 
\left(  \nu^3 (u_\nu''(t))^2+\frac{1}{\nu 
\phi(1/\nu)} \phi(|u_\nu'(t)|)\right) \hausmoins
\left(\{x \in A: d(x)=t\}\right) \ dt\\
&  &+\nu^3 \int_{0}^{\eta \nu} 
\left(|u_\nu''(t)||u_\nu'(t)| +
\frac{1}{2} 
(u_\nu'(t))^2 \vert\vert \nabla^2 d 
\vert\vert_\infty^2\right)\hausmoins\left(\{x \in A: d(x)=t\}\right) dt.
\end{eqnarray*}
We now claim that 
$$
F_\nu(\widetilde u_\nu, A) \leq
 \sup_{t \in (0,\eta \nu)} \left[F_\nu(u_\nu, (0, \eta \nu))\  \hausmoins(\{
x \in A:d(x)=t\}) \right]+ o(\nu).
$$
Indeed, since $\displaystyle \lim_{\nu \to 0^+} 
F_\nu(u_\nu, (-\gd,\gd))<+\infty$, we have
\[
\sup_{\nu} \nu^3 \int_{-\gd}^{\gd} (u_\nu'')^2 \ dt <+\infty,
\]
and by definition of $u_\nu$,
\[
\int_{0}^{\eta \nu} (u_\nu')^2~ dt 
=\frac{1}{\nu} \int_0^\eta (\psi')^2~ dt.
\]
Thus by Cauchy-Schwarz's inequality, 
\begin{align*}
\nu^3 
 \int_{0}^{\eta \nu}  
|u_\nu''||u_\nu'|~dt\le& 
\nu 
\left(\int_{0}^{\eta \nu} \nu^3  (u_\nu'')^2~dt\right)^{\frac{1}{2}}
 \left( \int_0^\eta (\psi')^2~ dt\right)^{\frac{1}{2}} = O(\nu).
\end{align*}
We can similarly bound the term
$\nu^3 \int_0^{\eta \nu} (u_\nu')^2~dt~ \vert\vert \nabla^2 d \vert\vert_\infty^2$. 
Finally since $\partial E $ is smooth,
\[\lim_{t \to 0^+}  \hausmoins(\{x \in A :d(x)=t\})
= \hausmoins(A \cap \partial E).
\]
Therefore, using the one-dimensional result in Section \ref{secequigamma}, we obtain
\[\limsup_{\nu \to 0^+} F_\nu
(\widetilde u_\nu,A) \le \sigma_a |s|^{\frac{2+a}{4-a}} \hausmoins(A \cap \partial E) +b.
\]
Letting $b \to 0^+$ we obtain the thesis.
\end{proof}

We can now prove the $\Gamma$-limsup, following the proof of \cite[Prop 3.6]{AlGe:01}.

\begin{Theorem}[{\bf $\Gamma$-limsup}]\label{pro:gsup}
 We have 
\[
\Gamma(L^2(A))-\limsup_{\nu \to 0^+} F_\nu(\cdot,A) \le \cF(\cdot,A),
\qquad A \in \mathcal A(\Omega).
\]
\end{Theorem}
\begin{proof}
 Thanks to the density Lemma \ref{teoantonin}, it is enough 
to prove the $\Gamma$-limsup
inequality  on those functions having a jump set 
contained in a finite union of facets of polytopes. Accordingly, we 
take
\[u= \sum_{i=1}^k s_i 1_{C_i},\]
where $k \in \mathbb N$, $s_i\in \R$
and  $C_i \subseteq A$ are closed polytopes with pairwise empty intersection of their
 interior parts. We can further assume that all $C_i$ are convex.  Set
\[C_{k+1}:= A \backslash \bigcup_{i=1}^k C_i\]
and 
\[\cI:=\{ (i,j) : i, j \in \{i,\dots, k+1\}, i<j, \hausmoins( C_i \cap C_j) \neq 0\}.\]
For $\gd>0$ let 
\[V^\gd:= \left\{ 
x \in A : {\rm dist}\Big(x, \bigcup_{i=1}^k \partial C_i\Big) < \gd 
\right\}\]
and 
\[V_{i j}^\gd := 
\big\{ x \in A : {\rm dist}(x, \partial C_i \cap \partial C_j) <\gd
\big\}, \qquad
(i,j) \in \cI.
\]
For any $(i,j)\in \cI$ choose 
functions $g_{ij} \in \cC^\infty_c(\Rn)$ such that
\[g_{ij}=1 \textrm{ on } V_{ij}^\gd, \qquad g_{ij}=0 \textrm{ on }
A \setminus
 V_{ij}^{2\gd},\]
and
\[
 \vert\vert  \nabla g_{ij}\vert\vert_{\infty} \le \frac{c}{\gd}, \qquad \vert\vert \nabla^2 g_{ij} \vert\vert_{\infty} \le \frac{c}{\gd^2},\]
where $c>0$ is a constant
 independent of $\gd$, $i$ and $j$. We now let
\[h_{ij}:= g_{ij} \left( \sum_{(l,m) \in \cI} g_{lm} + \prod_{(l,m)\in \cI} (1-g_{lm})\right)^{-1},\]
and 
\[h_0:= \prod_{(l,m) \in \cI} (1-g_{lm}) \left( \ \sum_{(l,m) \in \cI} g_{lm} + \prod_{(l,m)\in \cI} (1-g_{lm})\right)^{-1}.\]
These functions are defined so that
$$
h_0 + \sum_{(i,j)\in \cI} h_{ij}
 =1
$$
and
\[h_{ij}=g_{ij} \textrm{ on } V_{ij}^{2\gd} \backslash \bigcup_{(l,m) \neq (i,j)} V^{2\gd}_{lm} \qquad \textrm{and} \qquad h_0=1 \textrm{ on } A \backslash V^{2\gd}.\]

It can be verified that
\begin{align*}
 \vert\vert
\nabla h_{ij}\vert\vert_\infty\le 
\frac{c}{\delta}, \qquad & \vert\vert\nabla h_0
\vert\vert_\infty\le \frac{c}{\delta},\\[8pt]
\vert\vert\nabla^2 h_{ij}\vert\vert_\infty\le 
\frac{c}{\delta^2}, \qquad & \vert\vert\nabla^2 h_0\vert\vert_\infty\le \frac{c}{\delta^2}.
\end{align*}

Let $(u_{ij}^\nu) \subset H^2(\Omega)$ be the recovery sequence 
constructed as in Proposition \ref{gammasupE} related to 
\[u_{ij}:=\left\{\begin{array}{ll} s_i-s_j &\textrm{ on } S^{+}_{ij}\\
\\           0 &\textrm{ on } S^-_{ij},
          \end{array} \right.\]
where $S_{ij}$ 
is the hyperplane containing $C_i \cap C_j$ and 
\[S_{ij}^{\pm} := \{ x \in \Omega :  x=y\pm t\nu_{ij}, \ y \in S_{ij}, \ t \in \R^+\},\]
with $\nu_{ij}$ the unit internal normal to $C_i$ on $C_i \cap C_j$.
We fix then $\gd= \eta \nu$ where $\eta$ is given by the 
one-dimensional profile. Let also
\[u_\nu :=h_0 u
+ \sum_{(i,j) \in \cI} h_{ij} (u_{ij}^\nu +s_j).\]
Then $u_\nu \in H^2(\Omega)$ and $u_\nu \to u$ in $\Ldeu$
as $\nu \to 0^+$. 
Now on $V^{2\gd}_{ij} \backslash \bigcup_{(l,m) \neq (i,j)} 
V_{lm}^{2\gd}$, $u_\nu =(1-g_{ij}) u
+g_{ij}(u_{ij}^\nu+s_j)$ and since 
on $V^{2\gd}_{ij} \backslash V^{\gd}_{ij}$, $u_{ij}^\nu=u_{ij}$, 
we then have on this set $u_\nu=u$ and $\nabla u_\nu =0$. Thus

\begin{equation}\label{FnuVdelta}
F_\nu\left(u_\nu,  V^{2\gd}_{ij} \backslash \bigcup_{(l,m) 
\neq (i,j)} V_{lm}^{2\gd}\right)
=F_\nu\left(u_\nu,  V^{\gd}_{ij} \backslash \bigcup_{(l,m) 
\neq (i,j)} V_{lm}^{2\gd}\right)\le F_\nu\left(u_{ij}^\nu,
A \cap  V^{\gd}_{ij}\right).
\end{equation}
On $V^{2\gd}_{ij} \cap V^{2\gd}_{lm}$, we have (remember that 
$|\nabla u_{ij}^{\nu}|\le \vert\vert\psi'\vert\vert_\infty/\nu$)
\[|\nabla u_\nu| \le 
|u| \vert\nabla h_0\vert
+ \sum_{ (i,j) \in \cI}  
\left(|\nabla h_{ij}| |s_j| +|\nabla u^\nu_{ij}| \right)
\le \frac{C}{\nu},
\]
for a suitable constant $C>0$ independent of $\nu$.
Since $\cH^{\n-2}(C_i \cap C_j \cap C_l \cap C_m) < +\infty$
(where $\mathcal H^{n-2}$ is the $(n-2)$-dimensional
Hausdorff measure) we have $|V_{ij}^{2\gd}\cap V^{2\gd}_{lm}| = O(\nu^2)$ and thus 
\[
\int_{V_{ij}^{2\gd}\cap V^{2\gd}_{lm}} 
\  \frac{1}{\nu 
\phi(1/\nu)} \phi(|\nabla u_\nu|)~dx \le C \nu \frac{\phi(C/\nu)}{\phi(1/\nu)}.
\]
As $\lim_{\nu \to 0^+}  
\frac{\phi(C/\nu)}{\phi(1/\nu)}= C^a$, 
this shows that 
\[
 \lim_{\nu \to 0^+} \int_{V_{ij}^{2\gd}
\cap V^{2\gd}_{lm}} \ \frac{1}{\nu \phi(1/\nu)} \phi(|\nabla u_\nu|)~dx =0.
\]
Similarly,
\[
 \lim_{\nu \to 0^+} \int_{V_{ij}^{2\gd}\cap V^{2\gd}_{lm}} \ \nu^3 
|\nabla^2 u_\nu|^2~dx =0
\]
and therefore
\begin{equation}\label{FnuVinter}
 \lim_{\nu \to 0^+} F_\nu(u_\nu, V_{ij}^{2\gd}\cap V^{2\gd}_{lm})=0.
\end{equation}
Eventually on $A\backslash V^{2\gd}$, $u_\nu =u$ and $\nabla u =0$, thus 
\begin{equation}\label{FnuVcomp}
 F_\nu(u_\nu, A\backslash V^{2\gd})=0.
\end{equation}
Putting together 
\eqref{FnuVdelta}, \eqref{FnuVinter} and \eqref{FnuVcomp} 
we find by Proposition \ref{gammasupE},
\begin{align*}
 \lim_{\nu \to 0^+} F_\nu( u_\nu, A)\le &\sum_{(i,j) \in \cI} 
\lim_{\nu \to 0^+} F_\nu\left( u_\nu,  V^{2\gd}_{ij} \backslash \bigcup_{(l,m) \neq (i,j)} V_{lm}^{2\gd}\right)\\
  \le&\sum_{(i,j) \in \cI} 
\lim_{\nu \to 0^+} F_\nu( u^{ij}_\nu,  V_{ij}^\gd)\\
=&\sum_{(i,j) \in \cI} 
\lim_{\nu \to 0^+} F_\nu( u^{ij}_\nu,  A)\\
=&\sigma_a \sum_{(i,j) \in \cI} |s_i-s_j|^{\frac{2+a}{4-a}} 
\  \hausmoins(A \cap C_i \cap C_j)\\
=& \cF(u,A). 
\end{align*}
\end{proof}

\begin{Remark}\label{rem:final}\rm 
Differently with respect to 
the results of \cite{AlGe:01,Mo:01} 
(which rely on the density theorem of \cite{CoTo:})
using Lemma \ref{teoantonin} 
we obtain a full $\Gamma$-convergence result. Moreover, we have compactness in $BV$ for the $n$-dimensional problem which is a missing feature in 
most of the papers tackling similar problems such as \cite{AlGe:01,Mo:01}.
\end{Remark}

We conclude by mentioning that our theorems
are a first step toward a general approximation 
theorem in $SBV$ (see 
the related papers  \cite{Co:,CoTo:,AlBrGe:01,AlGe:01}) for 
free discontinuity functionals
having, in general, also a
bulk term.

\section*{Acknowledgments}
The first two authors 
 thank the hospitality of the Mathematisches Forschungsinstitut Oberwolfach,
which gave them the opportunity to start this research.

\end{document}